\documentclass{amsart}

\usepackage{epic}
\usepackage{eepic}

\newcommand     {\comment}[1]   {}
\newcommand{\mute}[2] {}
\newcommand     {\printname}[1] {}
\newcommand{\labell}[1] {\label{#1}\printname{#1}}

\def	\ss	{{\scriptsize}}
\def	\Z	{{\mathbb Z}}
\def	\ZZ	{{\mathbb Z}}
\def	\R	{{\mathbb R}}
\def	\RR	{{\mathbb R}}
\def	\C	{{\mathbb C}}
\def	\CP	{{\mathbb C}{\mathbb P}}
\def	\tCP	{{\widetilde{\CP^2}}}
\def	\tomega	{{\widetilde{\omega}}}
\def	\SO	{\operatorname{SO}}
\def	\GL	{\operatorname{GL}}
\def	\AGL	{\operatorname{AGL}}
\def	\Ham 	{\operatorname{Ham}}

\def    \half   {{\frac{1}{2}}}

\newtheorem {Theorem} {Theorem}
\newtheorem {Lemma} {Lemma}
\theoremstyle{remark}
\newtheorem*{Remark}{Remark}

\begin{document}

\title{Maximal tori in the symplectomorphism groups 
       of Hirzebruch surfaces}
\author{Yael Karshon}
\email{karshon@math.huji.ac.il}
\address{Institute of Mathematics, The Hebrew University of Jerusalem,
Giv'at Ram, Jerusalem 91904, Israel}

\thanks{This research was partially supported by the Israel Science
Foundation founded by The Academy of Sciences and Humanities.}


\begin{abstract}
We count the conjugacy classes of maximal tori in the groups of
symplectomorphisms of $S^2 \times S^2$ and of the blow-up of $\CP^2$ 
at a point.
\end{abstract}

\maketitle

Consider the group $\Ham (M,\omega)$ of Hamiltonian symplectomorphisms 
of a symplectic manifold.\footnote
{
Recall, a Hamiltonian symplectomorphism is one which can be 
connected to the identity by a path $\psi_t$ such that 
$\frac{d}{dt}\psi_t = X_t \circ \psi_t$ and $\iota(X_t)\omega = dH_t$
for a smooth $H_t \colon M \to \RR$.
}
A \emph{$k$-dimensional torus} in $\Ham(M,\omega)$ is a subgroup 
which is isomorphic to $(S^1)^k$. 
A \emph{maximal torus} is one which is not contained in any strictly larger 
torus. 

An action of $(S^1)^k$ on $(M,\omega)$ is called \emph{Hamiltonian}
if it admits a moment map, i.e., a map $\Phi \colon M \to \R^k$ 
such that $d\Phi_i = - \iota(\xi_i) \omega$ for all $i=1,\ldots,k$, 
where $\xi_1 , \ldots, \xi_k$ are the vector fields on $M$ 
that generate the action.

A Hamiltonian $(S^1)^k$-action defines a homomorphism 
from $(S^1)^k$ to $\Ham(M,\omega)$.
The action is effective if and only if this homomorphism is one to one.  
Its image is then a $k$-dimensional torus in $\Ham(M,\omega)$.
Every $k$-dimensional torus in $\Ham(M,\omega)$ is obtained in this way, 
and two Hamiltonian actions give the same torus if and only if 
they differ by a reparametrization of $(S^1)^k$. 
Tori in $\Ham(M,\omega)$ have dimension at most $\half \dim M$.
A Hamiltonian action of a $(\half \dim M)$-dimensional torus 
is called \emph{toric}.

\begin{Theorem} \labell{extends}
Let $(M,\omega)$ be a compact symplectic four-manifold.
Suppose that $\dim H^2(M,\R) \leq 3$ and $\dim H^1(M,\R) = 0$.
Then every Hamiltonian circle action on $(M,\omega)$
extends to a toric action.
\end{Theorem}

\begin{Remark}
Many symplectic four-manifolds do not admit Hamiltonian circle actions.
The theorem does not say anything about such manifolds.
\end{Remark}

Two tori, $T_1$ and $T_2$, in $\Ham(M,\omega)$ are \emph{conjugate} if 
there exists an element $g \in \Ham(M,\omega)$ such that $g T_1 g^{-1} = T_2$.
Two torus actions, viewed as homomorphisms $(S^1)^k \to \Ham(M,\omega)$,
give conjugate tori in $\Ham(M,\omega)$ 
if and only if they differ by an equivariant symplectomorphism 
composed with a reparametrization of $(S^1)^k$. 

On $S^2 \times S^2$, let $\omega_1$ and $\omega_2$ be the pullbacks
of the standard area form on $S^2$ via the two projection maps,
and let $\omega_{a,b} = a \omega_1 + b \omega_2$.
The standard $(S^1)^2$-action has the moment map image 
shown in Figure \ref{fig:min} on the left.

Let $\tCP$ be the blow-up of $\CP^2$ at a point. In it, let $E$ be
the exceptional divisor and let $L$ be a $\CP^1$ which is disjoint from $E$.
For $l > e > 0$, let $\tomega_{l,e}$ be a symplectic form such that
the symplectic area of $L$ is $l$ and the symplectic area of $E$ is $e$.
(We can construct $(\tCP,\tomega_{l,e})$ explicitly in several ways.
By \cite{McD}, it is unique up to symplectomorphism.)
The standard $(S^1)^2$-action (induced from the action on $\CP^2$)
has the moment map image shown in Figure \ref{fig:min} on the right.

\begin{figure}
\setlength{\unitlength}{0.00083333in}
\begingroup\makeatletter\ifx\SetFigFont\undefined%
\gdef\SetFigFont#1#2#3#4#5{%
  \reset@font\fontsize{#1}{#2pt}%
  \fontfamily{#3}\fontseries{#4}\fontshape{#5}%
  \selectfont}%
\fi\endgroup%
{\renewcommand{\dashlinestretch}{30}
\begin{picture}(3687,845)(0,0)
%
% left:
%
\path(75,150)(75,750)(1275,750) (1275,150)(75,150)
\put(-20,450){\ss $b$}
\put(450,0){\ss $a$}
%
% right:
%
\path(2175,150) (2175,750) (3075,750) (3975,150) (2175,150)
%\put(2070,450){\ss $b$}
%\put(2250,820){\ss $e=a-\frac{b}{2}$}
%\put(2250,-20){\ss $l=a+\frac{b}{2}$}
\put(2650,800){\ss $e$}
\put(2650,0){\ss $l$}
\put(3550,450){\ss slope=$-1$}
\end{picture} }
\caption{Moment map images for standard torus actions 
on $(S^2 \times S^2,\omega_{a,b})$ and $(\tCP,\tomega_{l,e})$.}
\labell{fig:min}
\end{figure}
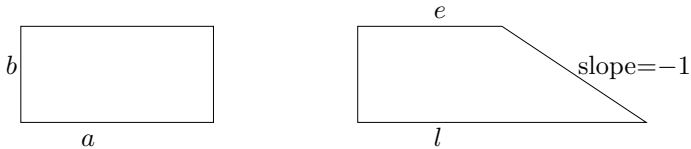

\begin{Theorem} \labell{max}
For $a \geq b > 0$, the number of conjugacy classes of maximal tori 
in $\Ham(S^2 \times S^2,\omega_{a,b})$ is 
%the number of integers $k$ that satisfy $0 \leq k < a/b$.
$\lceil a/b \rceil$.

For $l > e > 0$, the number of conjugacy classes of maximal tori 
in $\Ham(\tCP,\tomega_{l,e})$ is 
%the number of integers $k$ that satisfy $0 \leq k < \frac{e}{l-e}$.
$\lceil \frac{e}{l-e} \rceil$.
\end{Theorem}

\begin{Remark}
Here, $\lceil r \rceil$ denotes the smallest integer
greater than or equal to $r$.  We will prove Theorem \ref{max}
by enumerating the different conjugacy classes by the set 
of integers $k$ that satisfy $0 \leq k < r$
for $r = a/b$ and $r = \frac{e}{l-e}$, respectively.
\end{Remark}

\begin{Remark}
The topology of $\Ham(S^2 \times S^2,\omega_{a,b})$ also changes 
when the ratio $a/b$ crosses an integer.  
See \cite{gromov,abreu,abreu-McD,McD:ac,anjos}. 
Also see the remark at the end of the paper.
\end{Remark}

\begin{Remark}
There are infinitely many conjugacy classes of tori
in the group of contactomorphisms of an overtwisted $S^3$ 
or lens space. See \cite{lerman}.

For a count of the conjugacy classes of tori in the group of 
contactomorphisms of the pre-quantum line bundle over a 
Hirzebruch surface, see \cite{lerman2}.
\end{Remark}

We recall standard facts about Hamiltonian circle actions 
on compact symplectic manifolds:
By local arguments, each component of the fixed point set
is a submanifold of even dimension and has even index.
By Morse-Bott theory, each local extremum for the moment map 
is a global extremum, and these extrema are attained on connected sets.
See \cite[\S 32]{GS:book}. 
An \emph{interior fixed point} is a fixed point which is not a minimum
or maximum for the moment map.  A \emph{fixed surface} is 
a two dimensional connected component of the fixed point set.

\begin{Lemma} \labell{betti}
Let $M$ be a closed symplectic four-manifold with a Hamiltonian circle
action.  The dimension of $H^2(M)$ is equal to the number of interior
fixed points plus the number of fixed surfaces.  If $\dim H^1(M) = 0$,
each fixed surface has genus zero.
\end{Lemma}

\begin{proof}
We apply a standard Morse theory argument.
The moment map is a perfect Morse-Bott function whose critical points 
are the fixed points for the circle action \cite[\S 32]{GS:book}.
Therefore, $\dim H^j(M) = \sum \dim H^{j-i_F} (F)$,
where we sum over the connected components of the fixed point set, 
and where $i_F$ is the index of the component $F$.
The theorem follows by a simple computation of the summands, which is
summarized in the table below. 

\smallskip
\begin{tabular}{|c|c|c||c|c|c|c|c|}
\multicolumn{3}{c} \ &
\multicolumn{5}{c}{Contribution to } \\
\multicolumn{1}{c}{$F$} &
\multicolumn{1}{c}{$\Phi(F)$} &
\multicolumn{1}{c}{$i_F$} &
\multicolumn{1}{c}{$H^0(M)$} &
\multicolumn{1}{c}{$H^1(M)$} &
\multicolumn{1}{c}{$H^2(M)$} &
\multicolumn{1}{c}{$H^3(M)$} &
\multicolumn{1}{c}{$H^4(M)$} \\
\hline
fixed & minimal & $0$ & $1$ & $2 \operatorname{genus}(F)$ & $1$ & & \\
\cline{2-8}
surface & maximal & $2$ & & & $1$ & $2 \operatorname{genus}(F)$ & $1$ \\
\hline
isolated & minimal & $0$ & $1$ & & & & \\
\cline{2-8}
fixed & interior & $2$ & & & $1$ & & \\
\cline{2-8}
point & maximal & $4$ & & & & & $1$ \\
\hline
\end{tabular}
\end{proof}

\begin{proof}[Proof of Theorem \ref{extends}]
By \cite[Prop.~5.21]{karshon},
a Hamiltonian circle action extends to a toric action 
if and only if each fixed surface has genus zero
and each non-extremal level set for the moment map
contains at most two non-free orbits.

By Lemma \ref{betti}, since $\dim H^1(M)=0$, each fixed surface 
has genus zero.

By \cite[Theorem 5.1]{karshon}, if all the fixed points are isolated,
the circle action extends to a toric action.
Therefore, let us assume that there exists at least one fixed surface.
Suppose that the moment map attains its maximum on this surface;
the case of a minimum can be treated similarly.
By Lemma \ref{betti}, since $\dim H^2(M,\RR) \leq 3$, 
there exist at most two interior fixed points.

A \emph{$Z_k$-sphere} is a 2-sphere inside $M$ on which the circle 
acts by rotations with speed $k$.  A non-free orbit is either a fixed point
or belongs to a $Z_k$-sphere.  See \cite{audin} or \cite[Lemma 2.2]{karshon}.
A $Z_k$-sphere intersects each level set in at most one orbit. 
The north pole of a $Z_k$-sphere is an isolated, hence interior, fixed point
(because we assume that the maximal set of the moment map is not isolated).
Different $Z_k$-spheres have different north poles.
These considerations show that the number of non-free orbits
in an non-extremal level set for the moment map is at most the number
of interior fixed points. The theorem follows.
\end{proof}

For each non-negative integer $m$, consider the family of trapezoids
shown in Figure \ref{HirzAB}, parametrized by the height $b$ and
average width $a > \frac{m}{2} b$.
We call these \emph{standard Hirzebruch trapezoids}.
More generally, consider their images under the group $\AGL(2,\Z)$ 
of transformations of $\R^2$ of the form $x \mapsto Rx+v$ 
with $R \in \GL(2,\Z)$ and $v \in \R^2$.
These we call \emph{Hirzebruch trapezoids}.
Hirzebruch trapezoids modulo $\AGL(2,\Z)$ are in natural one-to-one
bijection with the set of parameters $(a,b,m)$
with $m$ a non-negative integer, $a$ and $b$ positive real numbers,
and $a > \frac{m}{2} b$.

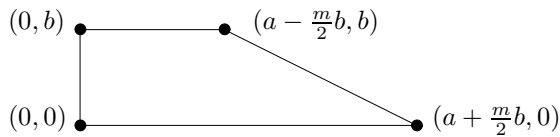
\begin{figure} 
\setlength{\unitlength}{0.00083333in}
\begingroup\makeatletter\ifx\SetFigFont\undefined%
\gdef\SetFigFont#1#2#3#4#5{%
  \reset@font\fontsize{#1}{#2pt}%
  \fontfamily{#3}\fontseries{#4}\fontshape{#5}%
  \selectfont}%
\fi\endgroup%
{\renewcommand{\dashlinestretch}{30}
\begin{picture}(3499,878)(0,-10)
\put(525,683){\blacken\ellipse{70}{70}}
\put(525,83){\blacken\ellipse{70}{70}}
\put(1425,683){\blacken\ellipse{70}{70}}
\put(2625,83){\blacken\ellipse{70}{70}}
\path(525,683)(1425,683)(2625,83) (525,83)(525,683)
\put(75,83){$\ss (0,0)$}
\put(75,683){$\ss (0,b)$}
\put(1600,683){$\ss (a-\frac{m}{2} b , b)$}
\put(2720,83){$ (a+\frac{m}{2}b , 0)$}
\end{picture}
}
\caption{A standard Hirzebruch trapezoid}
\labell{HirzAB}
\end{figure}

Take a polygon in $\R^2$. Let $u_1,\ldots,u_m$ be normal vectors
to its edges, in counterclockwise order, pointing inward. 
The polygon is a \emph{Delzant polygon} if we can choose these normals
so that $u_i \in \Z^2$ and $\det(u_i u_{i+1}) =1$ for all $i$,
the indices taken cyclically.\footnote{
More generally, the Delzant condition for a polytope in $\R^n$
is that exactly $n$ facets meet at every vertex
and the normals to these facets can be chosen to be generators of $\Z^n$.
}

\begin{Lemma} \labell{all are hirz}
Every Delzant polygon with four edges is a Hirzebruch trapezoid.
\end{Lemma}

\begin{proof}
Fix a Delzant polygon with four edges. 
Let $u_1$, $u_2$, $u_3$, and $u_4$ be normals to its edges
that satisfy the above Delzant condition.
Because $\det(u_1 u_2) = 1$, we may assume, without loss of generality, 
that $u_1 = (1,0)$ and $u_2 = (0,1)$.
The conditions $\det (u_2 u_3) = 1$ and $\det (u_4 u_1) = 1$ imply
that $u_3 = (-1,k)$ and $u_4 = (l,-1)$ for some $k,l \in \Z$.
The condition $\det (u_3 u_4) = 1$ then implies $kl = 0$.
Each of the cases $k=0$ and $l=0$ gives a Hirzebruch trapezoid,
as shown in Figure \ref{two Hirz}.
\end{proof}

\begin{figure}
\setlength{\unitlength}{0.00063333in}
\begingroup\makeatletter\ifx\SetFigFont\undefined%
\gdef\SetFigFont#1#2#3#4#5{%
  \reset@font\fontsize{#1}{#2pt}%
  \fontfamily{#3}\fontseries{#4}\fontshape{#5}%
  \selectfont}%
\fi\endgroup%
{\renewcommand{\dashlinestretch}{30}
\begin{picture}(4812,2739)(-150,-10)
%
% left:
%
\path (150,12) (150,912) (1050,2712) (1050,12) (150,12)
\path(150,462)(375,462)
\path(255.000,432.000)(375.000,462.000)(255.000,492.000)
\path(450,1512)(600,1437)
\path(479.252,1463.833)(600.000,1437.000)(506.085,1517.498)
\path(1050,912)(825,912)
\path(945.000,942.000)(825.000,912.000)(945.000,882.000)
\path(600,12)(600,237)
\path(630.000,117.000)(600.000,237.000)(570.000,117.000)
% text on left:
\put(0,1812){\ss slope}
\put(50,1572){\ss $=l$}
\put(225,537){\ss $u_1$}
\put(675,87){\ss $u_2$}
\put(825,737){\ss $u_3$}
\put(525,1312){\ss $u_4$}
%
% right:
%
\path (1950,12) (1950,912) (4650,912) (2850,12) (1950,12)
\path(1950,462)(2175,462)
\path(2055.000,432.000)(2175.000,462.000)(2055.000,492.000)
\path(2850,912)(2850,687)
\path(2820.000,807.000)(2850.000,687.000)(2880.000,807.000)
\path(3450,312)(3375,462)
\path(3455.498,368.085)(3375.000,462.000)(3401.833,341.252)
\path(2400,12)(2400,237)
\path(2430.000,117.000)(2400.000,237.000)(2370.000,117.000)
\put(3750,312){\ss slope $=\frac{1}{k}$}
\put(2025,537){\ss $u_1$}
\put(2475,87){\ss $u_2$}
\put(3450,462){\ss $u_3$}
\put(2925,687){\ss $u_4$}
\end{picture}
}
\caption{Hirzebruch trapezoids}
\labell{two Hirz}
\end{figure}
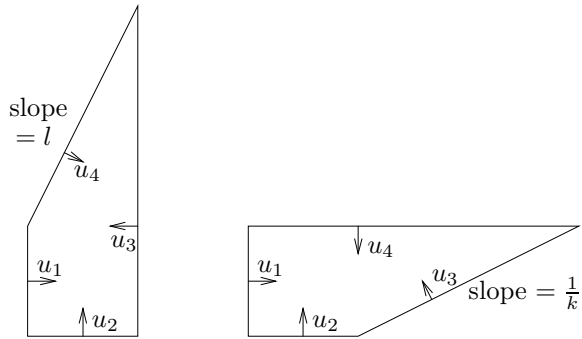

The convexity theorem of Atiyah, Guillemin and Sternberg, [At,Gu-St], 
states that, for a Hamiltonian torus action an 
a compact symplectic manifold, the image of the moment map 
is a convex polytope. 
By Delzant's classification of Hamiltonian toric actions [De],\footnote{
Also see [L-T].
}
the moment map images for \emph{toric} actions 
are precisely the Delzant polytopes,
and two toric actions are equivariantly symplectomorphic
if and only if their moment map images are translates of each other.

The symplectic four-manifolds that correspond to Hirzebruch trapezoids 
are \emph{Hirzebruch surfaces}. 
See \cite{audin} or \cite[section 6.3]{karshon}.
Specifically, when $m=0$ or $m=1$, these are
$(S^2 \times S^2,\omega_{a,b})$ and $(\tCP,\tomega_{l,e})$, respectively.

\begin{Lemma} \labell{hirz even}
A Hirzebruch surface which corresponds to an integer $n \geq 2$
is symplectomorphic to the Hirzebruch surface which corresponds
to the integer $n-2$ with the same parameters $a$ and $b$.
\end{Lemma}

\begin{proof}
Let $m=n-1$.
Take two Hirzebruch trapezoids with the same parameters $a$ and $b$,
corresponding to the integers $m-1$ and $m+1$. After transforming by 
appropriate elements of $\AGL(2,\Z)$, they can be brought to the 
forms shown in Figure \ref{m pm 1}.

In \cite[\S 2.1]{karshon} we associate a labeled graph 
to every Hamiltonian $S^1$-space, such that two spaces are isomorphic
if and only if their graphs are isomorphic \cite[Theorem 4.1]{karshon}. 
If the $S^1$-action is obtained by restriction of a toric action, 
the graph can be easily read from the Delzant polygon.
See \cite[\S 2.2]{karshon}. The two polygons in Figure~\ref{m pm 1} 
give rise to the same graph (shown in Figure \ref{m pm 1} on the right).
Therefore, the spaces are $S^1$-equivariantly symplectomorphic.
\end{proof}

\begin{Remark} 
Lemma \ref{hirz even} (which I learned from S.\ Tolman) seems to be 
well known.  
For instance, the new circle action on $S^2 \times S^2$
obtained by identifying this space with the Hirzebruch surface
with parameters $(a,b,m=2)$ plays a role in \cite{McD:examples}.
\end{Remark}

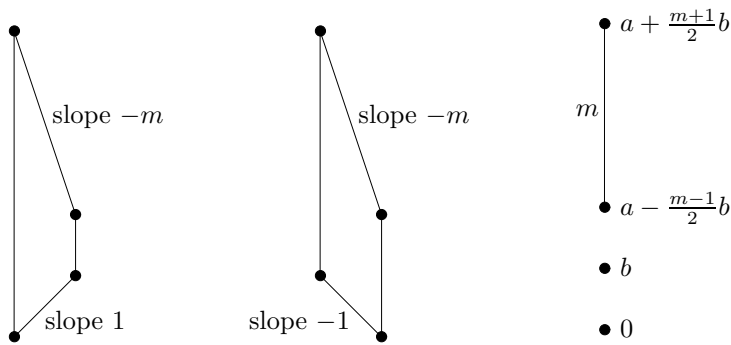
\begin{figure} 
\setlength{\unitlength}{0.00053333in}
\begingroup\makeatletter\ifx\SetFigFont\undefined%
\gdef\SetFigFont#1#2#3#4#5{%
  \reset@font\fontsize{#1}{#2pt}%
  \fontfamily{#3}\fontseries{#4}\fontshape{#5}%
  \selectfont}%
\fi\endgroup%

\ \hfill
\begin{picture}(4156,3039)(0,-10)
\path(12,3012)(612,1212)(612,612) (12,12)(12,3012)
\put(12,3012){\blacken\ellipse{100}{100}}
\put(612,1212){\blacken\ellipse{100}{100}}
\put(612,612){\blacken\ellipse{100}{100}}
\put(12,12){\blacken\ellipse{100}{100}}
\path(3012,3012)(3612,1212)(3612,12)(3012,612)(3012,3012)
\put(3012,3012){\blacken\ellipse{100}{100}}
\put(3612,1212){\blacken\ellipse{100}{100}}
\put(3612,12){\blacken\ellipse{100}{100}}
\put(3012,612){\blacken\ellipse{100}{100}}
\put(387,2112){\ss slope $-m$}
\put(3387,2112){\ss slope $-m$}
\put(312,87){\ss slope $1$}
\put(2312,87){\ss slope $-1$}
\end{picture}
\hfill
\begin{picture}(494,3181)(0,-10)
\put(83,3083){\blacken\ellipse{100}{100}}
\put(83,1283){\blacken\ellipse{100}{100}}
\put(83,683){\blacken\ellipse{100}{100}}
\put(83,83){\blacken\ellipse{100}{100}}
\path(83,3083)(83,1283)
\put(233,8){\ss $0$}
\put(233,608){\ss $b$}
\put(233,3008) {\ss $a + \frac{m+1}{2} b$}
\put(233,1208) {\ss $a - \frac{m-1}{2} b$}
\put(-200,2183){\ss $m$}
\end{picture}
\hfill \ 
\caption{Hirzebruch trapezoids with integers $m-1$ and $m+1$,
width $b$ and average height $a$.}
\labell{m pm 1}
\end{figure}

\begin{Lemma} \labell{distinct}
Among the symplectic manifolds 
$(S^2 \times S^2 , \omega_{a,b})$, for $a \geq b > 0$,
and $(\tCP, \tomega_{l,e})$, for $l > e > 0$,
no two are symplectomorphic.
\end{Lemma}

\begin{proof}
The manifolds $S^2 \times S^2$ and $\tCP$ are not homeomorphic. 
For instance, the self intersection of the exceptional divisor 
in $\tCP$ is $-1$ whereas in $S^2 \times S^2$ every class in $H_2$ 
has an even self intersection.

A diffeomorphism of $S^2 \times S^2$ 
acts on $H^2(S^2 \times S^2) = \Z^2$ 
by a $2 \times 2$ matrix of integers, with determinant $\pm 1$,
which preserves the intersection form
$ \left[ \begin{array} {cc}  0 & 1 \\ 1 & 0 \end{array} \right]. $
There are exactly sixteen matrices with this property; they are
$$ \left[ \begin{array} {cc}  \pm 1 & 0 \\ 0 & \pm 1 \end{array} \right]
\quad \text{and} \quad
 \left[ \begin{array} {cc}  0 & \pm 1 \\ \pm 1 & 0 \end{array} \right].$$
These cannot take $\omega_{a,b}$ to a different $\omega_{a',b'}$
with $a \geq b > 0$ and $a' \geq b' > 0$.

A diffeomorphism of $\tCP$ 
acts on $H_2(\tCP) = \ZZ L \oplus \ZZ E \cong \ZZ^2$
by a $2 \times 2$ matrix of integers, with determinant $\pm 1$,
which preserves the intersection form 
$ \left[ \begin{array}{rr} 1 & 0 \\ 0 & -1 \end{array} \right]$.
There are exactly four matrices with this property; they are
$\left[ \begin{array}{rr} \pm 1 & 0 \\ 0 & \pm 1 \end{array} \right]$.
These cannot take $\tomega_{l,e}$ to a different $\tomega_{l',e'}$
with $l > e > 0$ and $l' > e' > 0$.
\end{proof}

\begin{proof}[Proof of Theorem \ref{max}]
By Theorem \ref{extends}, we only need to consider two dimensional tori.

The conjugacy classes of $2$-tori in $\Ham(M,\omega)$
for \emph{all} possible symplectic four-manifolds $(M,\omega)$
are given by all the Delzant polygons in $\R^2$,
modulo $\AGL(2,\Z)$-congruence.
This follows immediately from Delzant's classification 
of Hamiltonian toric actions [De] and the fact that
a reparametrization of the 2-torus $S^1\times S^1$ 
transforms the moment map image by an element of $\GL(2,\Z)$.

Hence, to find the conjugacy classes of 2-tori in the symplectomorphism
group of a particular symplectic four-manifold $(M,\omega)$,
we must identify which Delzant polygons correspond to a space 
which is symplectomorphic to $(M,\omega)$, and we must take 
these polygons modulo $\AGL(2,\Z)$-congruence.

The number of edges of a Delzant polygon is equal to $2$ 
plus the second Betti number of the corresponding space. 
See Lemma \ref{betti} and \cite[section 2.2]{karshon}.
Therefore, when $M = S^2 \times S^2$ or $M = \tCP$,
we only need to consider Delzant polygons with four edges.
By Lemma \ref{all are hirz}, we only need to consider Hirzebruch trapezoids.

Consider the Hirzebruch trapezoid with parameters $m \geq 0$
and $a \geq b > 0$.  Iterate Lemma \ref{hirz even}. 
If $m$ is even, the corresponding Hirzebruch surface 
is symplectomorphic to $(S^2 \times S^2,\omega_{a,b})$.
If $m$ is odd, it is symplectomorphic to $(\tCP,\tomega_{l,e})$,
where $l = a + \frac{b}{2}$ and $e = a - \frac{b}{2}$,
so that the corresponding trapezoid 
(shown in Figure \ref{fig:min} on the right)
still has height $b$ and average width $a$. 

Numbers $a \geq b$ can occur as the average width and the height
of a Hirzebruch trapezoid with integer parameter $m \geq 0$
if and only if $\frac{a}{b} > \frac{m}{2}$.  
When $m=2k$ is even, this becomes the condition $k<\frac{a}{b}$.
When $m=2k+1$ is odd, 
it becomes the condition $k<\frac{e}{l-e}$, with $l$ and $e$ as above.

We have found the required number of distinct torus actions on each 
of the symplectic manifolds $(S^2 \times S^2 , \omega_{a,b})$,
$a \geq b > 0$, and $(\tCP,\tomega_{l,e})$, $l > e > 0$.
By Lemma \ref{distinct}, this accounts for all possible toric actions
on each of these symplectic manifolds.
\end{proof}

\begin{Remark}
To each non-negative integer $m$ we have associated a torus action
on $S^2 \times S^2$ if $m$ is even and on $\widetilde{\CP^2}$
if $m$ is odd, for appropriate ranges of values of symplectic forms. 
Each of these torus actions is in fact obtained from an action of a 
larger, non-abelian, compact Lie group by restricting to its maximal
torus.  For instance, the standard actions of $\SO(3) \times \SO(3)$
on $S^2 \times S^2$ and of $U(2)$ on $\widetilde{\CP^2}$
restrict to the torus actions that correspond to the integers $m=0$
and $m=1$, respectively.
More generally, given a non-negative integer $m$, consider the quotient
of $S^3 \times S^2$ by the circle action
$\lambda \colon (z,p) \mapsto (z \lambda, \lambda^m \cdot p)$,
where $S^1$ is the circle group of complex numbers of norm one,
$\lambda \in S^1$ acts on $z \in S^3 \subset \C^2$ by scalar multiplication,
and $p \mapsto \lambda \cdot p$ is the circle action on
$S^2 \subset \R^3$ by rotations.
This space admits natural K\"ahler structures with which it is 
symplectomorphic (but not biholomorphic) to $S^2 \times S^2$ if $m$ is even
and to $\widetilde{\CP^2}$ if $m$ is odd, with appropriate symplectic forms.
On this space there is a natural action of the quotient 
of $U(2) \times S^1$ by the subgroup
$\{ (aI,a^m) \mid a \in S^1 \}$, where $I \in U(2)$ is the identity matrix.
Note that this quotient is a central extension of $\SO(3)$ by $S^1$.
The torus action that corresponds to the integer $m$ comes from this action.

In this context we recall that the only four dimensional compact symplectic 
$\SO(3)$ manifolds are $S^2 \times S^2$ and $\CP^2$, by \cite{iglesias}.

We also recall that the Hamiltonian symplectomorphism group 
of $S^2 \times S^2$ with equal areas of the two factors
retracts to $\SO(3) \times \SO(3)$ by \cite{gromov}
and that the aforementioned compact Lie subgroups 
of the Hamiltonian symplectomorphism groups of $S^2 \times S^2$ 
and of $\widetilde{\CP^2}$ in some sense carry an essential part 
of the topology of these symplectomorphism groups  
\cite{abreu,abreu-McD,McD:ac}.
\end{Remark}

%%%%%%%%%%%%%%%%%%%%%%%%%%%%%%%%%%%%%%%%%%%%

\subsection*{Acknowledgement}
I thank Sue Tolman for calling my attention, years ago, 
to the argument of Lemma \ref{hirz even}. 
I thank F.~Lalonde for calling my attention, also years ago,
to the question of finding different actions on the same symplectic manifold.
I thank an anonymous referee for helpful comments.
I thank Eugene Lerman for pushing me to publish this paper 
and for helpful comments.

\end{document}